\documentclass[12pt]{amsart}

\usepackage{amsmath,amsfonts,amssymb,amscd}
\begin{document}
\renewcommand{\thefootnote}{\fnsymbol{footnote}}
\pagestyle{plain}

%%%%%%%%%%%%%%%%%%%%%%%%
\title{Asymptotics of
 polybalanced metrics \\
under relative stability constraints}
\author{Toshiki Mabuchi${}^*$}
%\address{Department of Mathematics, Graduate School of Science, Osaka
%University,  Toyonaka, Osaka, 560-0043 Japan
%}
\maketitle
%$$
%\align
%&x=y
%\tag"(1.1)" \\
%&z=w \tag"(1.2)"
%\endalign
%$$
%%%%%%%%%%%%%%%%%%%%%%%%
\footnotetext{ ${}^{*}$Supported 
by JSPS Grant-in-Aid for Scientific Research (A) No. 20244005.}
%The result in this paper is announced in ``Extremal metrics and 
%K\"ahler-Ricci flow'' held in Pisa. }
%%%%%%%%%%%%%%%%%%%%%%%%
\abstract
Under the assumption of asymptotic relative Chow-stability for polarized algebraic manifolds $(M, L)$, 
a series of weight\-ed balanced metrics $\omega_m$, $m \gg 1$, called 
polybalanced metrics, are obtained from complete linear systems $|L^{m}|$ on $M$.
Then the asymptotic behavior of the weights 
as $m\to \infty$ will be studied. 
\endabstract
%%%%%%%%%%%%%%%%%%%%%%

\section{Introduction}

In this paper, we shall study relative Chow-stability (cf.~\cite{Ma}; see also \cite{Sz}) for polarized algebraic manifolds $(M,L)$ from the viewpoints of the  existence problem of extremal K\"ahler metrics.
As balanced metrics are obtained from Chow-stability on polarized algebraic manifolds,  
our relative Chow-stability similarly provides us with
a special type of weighted balanced metrics called {\it polybalanced metrics}.
As a crucial step in the program of \cite{M3}, we here study the asymptotic behavior of the weights for such polybalanced metrics. 

\medskip
By a {\it polarized algebraic manifold} $(M,L)$, we
 mean a pair of a connected 
projective algebraic manifold $M$ 
and a very ample holomorphic line bundle $L$
over $M$. 
For a maximal connected linear algebraic subgroup $G$ of 
the group $\operatorname{Aut}(M)$ of all holomorphic automorphisms of $M$, 
let $\frak g := \operatorname{Lie} G$ denote its Lie algebra.
%For a maximal connected linear algebraic subgroup $G$ of 
%the group $\operatorname{Aut}(M)$ of all holomorphic automorphisms of $M$, 
%the Chevalley decomposition allows us to write 
%$$
%G \; =\; R_{\Bbb C} \ltimes U
%$$
%as a semidirect product of a reductive algebraic group $R_{\Bbb C}$ 
%and  the unipotent radical $U$ of $G$. 
%Let  $\frak g :=\operatorname{Lie}G$ and 
%$\frak r := \operatorname{Lie}R_{\Bbb C}$ be the Lie algebras of 
%$G$ and $R_{\Bbb C}$, respectively. For the center $\frak z$ of $\frak r$,
%fixing an arbitrary maximal toral subalgebra ${\frak t}_{\operatorname{max}}$ of 
%$\frak r$,
%we see that ${\frak t}_{\operatorname{max}}$ is a maximal toral subalgebra of $\frak g$
%satisfying $\frak z \subset {\frak t}_{\operatorname{max}}$.
%Consider the series of the Lie algebra characters (cf.~\cite{F1}) 
%$$
%\mathcal{F}_p :\frak g \to \Bbb C, 
%\qquad p=1,2,\dots,n,
%$$
%defined as obstructions to asymptotic Chow semistability of $(M,L)$,
%where $\mathcal{F}_1$ is the classical Futaki character of $M$.
%Define a subspace $\frak a$ of $\frak z$, defined over $\Bbb Q$, consisting of all 
%$A\in \frak z$ such that
%$\mathcal{F}_p (A) =0$ for all\; $p = 1,2,\dots,n$. 
Since the infinitesimal $\frak g$-action on $M$ 
lifts to an infinitesimal bundle $\frak g$-action on $L$, by setting
$$
V_m:= H^0(M, L^{ m}), \qquad m =1,2,\dots,
$$
we view $\frak g$ as a Lie subalgebra of $\frak{sl}(V_m)$. 
We now define a symmetric bilinear form $ \langle \;, \,\rangle^{}_m$ 
on $\frak{sl}(V_m)$ by
$$
 \; \langle X, Y\rangle^{}_m \; =\; \operatorname{Tr}(XY)/m^{n+2},
 \qquad X, Y\in \frak{sl}(V_m), 
$$
where the asymptotic limit of $ \langle \;, \,\rangle^{}_m$ as $m\to \infty$ often 
plays an important role in the study of K-stability. In fact one can show
$$
\langle X, Y\rangle^{}_m \; =\; O(1)
$$
by using the equivariant Riemann-Roch formula, see \cite{D2} and \cite{Sz}.
%Since $ \langle \;, \,\rangle^{}_m$ restricted to 
%the Lie subalgebra $\frak z$ of $\frak{sl}(V_m)$
%is nondegenerate for each positive integer $m$, we can define a complex Lie algebra
%$$
%{\frak b}_{m} := 
%\frak a^{\perp m}
%$$
%as the orthogonal complement, defined over $\Bbb Q$,   
%of $\frak a$ in $\frak z$ consisting of all $B\in \frak z$ such that 
%$\langle A,B\rangle_m^{} = 0$ for all $A \in \frak a$.
%Let $\frak t_{\operatorname{min}}$ 
%denote the complex Lie subalgebra, defined over $\Bbb Q$,
%of $\frak z$ generated by all 
%$$
%\frak b_m, \qquad m = 1,2,\dots, 
%$$
%in the center $\frak z$.
Let $T$ be an algebraic torus in $\operatorname{SL}(V_m)$ 
such that the corresponding Lie algebra $\frak t := \operatorname{Lie} T$ 
satisfies 
$$
\frak t \; \subset\; \frak g.
$$
%$$
%\frak t_{\operatorname{min}}\subset  \frak t.
%$$
Then by the 
$T$-action on $V_m$, we can write the vector space $V_m$ as a direct sum 
of $\frak t$-eigenspaces:
$$
V_m \; =\; \bigoplus_{k=1}^{\nu_m}\; V(\chi^{}_k ),
$$
where $V(\chi_k ) :=
\{\,v\in V_m\,;\, g\cdot v = \chi^{}_k (g) v \text{ for all }g\in T\}$
for mutually distinct multiplicative characters
$\chi^{}_k \in \operatorname{Hom}(T, \Bbb C^* )$, 
$k = 1,2,\dots, \nu_m$.  

\medskip
To study $V_m$, let $\omega_m$ be a K\"ahler metric in the class $c_1(L)_{\Bbb R}$, 
and choose a Hermitian metric $h_m$ for $L$ such that 
$\omega_m = c_1(L; h_m)_{\Bbb R}$. We now endow $V_m$ 
with the Hermitian $L^2$ inner product 
on $V_m$ 
defined by
$$
(u, v)_{L^2}:=  \int_M (u,v)^{}_{h_m} \, \omega_m^{\,n}, 
\qquad u,v\in V_m,
\leqno{(1.1)}
$$
where $(u,v)_{h_m}$ denotes the pointwise Hermitian pairing of $u$, $v$ in terms of $h_m$. 
Then by this $L^2$ inner product, we have $V(\chi_k ) \perp V(\chi_{k'} )$, $k\neq k'$. 
Put $N_m := \dim V_m$ and $n_k:= \dim_{\Bbb C} V(\chi_k )$. For each $k$, by choosing an orthonormal basis 
$\{\,\sigma_{k,i}\,;\, i=1,2,\dots,n_k\,\}$ for $V(\chi_k )$, we put
$$
B_{m,k} (\omega_m ): = \; \sum_{i=1}^{n_k}\; |\sigma_{k,i}|_{h_m}^2,
$$ 
where $|u|^2_{h_m} := (u,u)_{h_m}$ for each $u\in V_m$.
Then $\omega_m$ is called a {\it polybalanced metric},  if there exist 
real constants $\gamma_{m,k} >0$ such that
$$
B^{\circ}_m (\omega_m ) = \; \sum_{k=1}^{\nu_m}\; \gamma_{m,k}\, 
B_{m,k} (\omega_m )
\leqno{(1.2)}
$$
is a constant function on $M$.
Here $\gamma_{m,k}$ are called the {\it  weights} of the polybalanced metric $\omega_m$.
On the other hand,
$$
B^{\bullet}_m (\omega_m )\,:=\; \sum_{k=1}^{\nu_m}\; 
B_{m,k} (\omega_m )
$$
is called the {\it $m$-th asymptotic Bergman kernel} of $\omega_m$.
A smooth real-valued function $f\in C^{\infty}(M)_{\Bbb R}$ on the K\"ahler manifold 
$(M,\omega_m)$ is said to be {\it Hamiltonian} if there exists a holomorphic vector field 
$X\in \frak g$ on $M$ such that $i_X \omega_m = \sqrt{-1}\,\bar{\partial}f$.
Put $N_m':= N_m/ c_1(L)^n[M]$.
In this paper, as the first step in \cite{M3}, we shall show the following:

\medskip\noindent
{\bf Theorem A:} 
{\em 
For a polarized algebraic manifold $(M,L)$ and an algebraic torus $T$ as above,
assume that $(M, L)$ is asymptotically Chow-stable relative to $T$.
Then for each $m \gg 1$, there exists a polybalanced metric $\omega_m$ in the class $c_1(L)_{\Bbb R}$ such that
 $\gamma_{m,k} = 1 + O(1/m)$, i.e.,
 $$
 |\gamma_{m,k} - 1| \;\leq \; C_1 /m, 
 \qquad k= 1,2,\dots, \nu_m; \;\; m \gg 1,
 \leqno{(1.3)}
 $$
 for some positive 
constant $C_1$ independent of $k$ and $m$. Moreover, there exist 
uniformly $C^0$-bounded functions $f_m \in C^{\infty}(M)^{}_{\Bbb R}$  on $M$ 
such that
$$
B^{\bullet}_m (\omega_m )\; =\; N'_m + f_m m^{n-1} + O(m^{n-2})
\leqno{(1.4)}
$$
and that each $f_m$ is a Hamiltonian function on $(M, \omega_m )$ satisfying
$i_{X_m} \omega_m = \sqrt{-1}\,\bar{\partial} f_m$
for some holomorphic 
vector field $X_m \in \frak t$ on $M$.
}

\medskip
In view of \cite{MN}, this theorem and
the result of Catlin-Lu-Tian-Yau-Zelditch (\cite{L},\cite{T1},\cite{Ze}) 
 allow us to obtain an approach (cf.~\cite{M3}) to an extremal K\"ahler version of Donaldson-Tian-Yau's conjecture.
On the other hand, as a corollary to Theorem A, we obtain the following:

\medskip\noindent
{\bf Corollary B:} 
{\em Under the same assumption as in Theorem A, suppose further that the 
classical Futaki character 
$\mathcal{F}_1: \frak g \to \Bbb C$ for $M$ vanishes on $\frak t$. Then for each $m \gg 1$, there exists a polybalanced metric $\omega_m$ in the class $c_1(L)_{\Bbb R}$ such that
 $\gamma_{m,k} = 1 + O(1/m^2)$. In particular
$$
B^{\bullet}_{m} (\omega_m )\; =\; N'_m +  O(m^{n-2}).
$$
}

\section{Asymptotic relative Chow-stability}

By the same notation as in the introduction,
we consider the algebraic subgroup $S_m$ 
of $\operatorname{SL}(V_m)$ defined by
$$
S_m \; :=\; \prod_{k=1}^{\nu_m}\; \operatorname{SL}(V (\chi_k )),
$$
where the action of each $\operatorname{SL}(V (\chi_k ))$ on $V_m$ 
fixes $V(\chi_i )$ if $i \neq k$.
Then the centralizer $H_m$ of $S_m$ in $\operatorname{SL}(V_m)$
consists of all diagonal matrices in $\operatorname{SL}(V_m)$
acting on each $V (\chi_k )$ by constant scalar multiplication.
Hence the centralizer $Z_m(T)$ of $T$ in $\operatorname{SL}(V_m)$ is $H_m\cdot S_m$ 
with Lie algebra
$$
\frak z_{m}(\frak t ) \; =\; \frak h_m \, +\, \frak s_m,
$$
where $\frak s_m := \operatorname{Lie} S_m$ and $\frak h_m := \operatorname{Lie}H_m$.
For the exponential map defined by
$\frak h_m \owns X \mapsto \exp (2\pi \sqrt{-1}\, X) \in H_m$, 
let $(\frak h_m)_{\Bbb Z}$ denote its kernel.
Regarding $(\frak h_m)_{\Bbb R}:=(\frak h_m)_{\Bbb Z}\otimes_{\Bbb Z}\Bbb R$ 
as a subspace of $\frak h_m$, we have a real structure on $\frak h_m$, 
i.e., an involution
$$
\frak h_m \owns X \mapsto \bar{X}\in \frak h_m
$$ 
defined as the associated complex conjugate of $\frak h_m$ 
fixing $(\frak h^{}_m)^{}_{\Bbb R}$. 
We then have a Hermitian metric $(\;,\,)_m$ on $\frak h_m$ 
by setting
$$
(X, Y)_m \; =\; \langle X, \bar{Y}\rangle^{}_m,  \qquad X, Y\in \frak h_m. 
\leqno{(2.1)}
$$
For the orthogonal complement $\frak t^{\perp}$ 
of $\frak t$ in $\frak h_m$ in terms of this Hermitian metric, let $T^{\perp}$ denote the corresponding algebraic torus 
in $H_m$. We now define an algebraic subgroup $G_{m}$ 
of $Z_m (T)$ by
$$
G_{m} :=  T^{\perp}\cdot S_m. 
\leqno{(2.2)}
$$
For the $T$-equivariant Kodaira embedding
$\Phi_m : M \hookrightarrow \Bbb P^*(V_m)$
associated to the complete linear system $| L^m |$ on $M$,
let $d(m)$ denote the degree of the image $\Phi_m (M)$ in the projective space 
$\Bbb P^*(V_m)$. 
For the dual space $W^*_m$ of $W_m:=  S^{d(m)} (V_m)^{\otimes n+1}$,
we have the Chow form 
$$
0 \,\neq\, \hat{M}_m\,\in\, W^*_m
$$ 
for the irreducible reduced algebraic cycle 
$\Phi_m (M)$ on $\Bbb P^* (V_m )$,  
so that the corresponding element $[\hat{M}_m ]$ in $\Bbb P^* (W_m )$
is the Chow point for the cycle $\Phi_m (M)$.
Consider the natural action of $\operatorname{SL}(V_m)$ on $W^*_m$ 
induced by the action of $\operatorname{SL}(V_m)$ on $V_m$.

\medskip\noindent
{\bf Definition 2.3.}
 (1) $(M,L^m)$ is said to be {\it Chow-stable relative to} $T$ if 
the orbit $G_m\cdot \hat{M}_m$ is closed in $W^*_m$.
 
 \smallskip\noindent
 (2) $(M,L)$ is said to be {\it  asymptotically Chow-stable relative to} $T$ if $(M, L^{m})$ 
 is Chow stable relative to $T$ for each integer $m \gg 1$.
 
 \section{Relative Chow-stability for each fixed $m$}
 
In this section, 
we consider a polarized algebraic manifold $(M,L)$
under the assumption that $(M,L^m)$ is Chow stable 
 relative to $T$ for a fixed positive integer $m$. 
Then we shall show that a polybalanced metric 
$\omega_m$ exists in the class $c_1(L)_{\Bbb R}$.
 
\medskip
The space
$\Lambda_m := \{\,\lambda = (\lambda_1, \lambda_2,\dots, \lambda_{\nu_m})\in \Bbb C^{\nu_m}\,;\,
\Sigma_{k=1}^{\nu_m}\, n_k\lambda_k =0\, \}$ and the Lie algebra $\frak h_m$ 
are identified by an isomorphism
$$
\Lambda_m \;  \cong \;  {\frak h}_m, 
\qquad \lambda \leftrightarrow X_{\lambda},
\leqno{(3.1)}
$$
with $(\frak h_m )_{\Bbb R}$ corresponding to the set $(\Lambda_m)_{\Bbb R}$ 
of the real points in $\Lambda_m$,
where $X_{\lambda}$ is the endomorphism of $V_m$ defined by
$$
X_\lambda \; := \; \bigoplus_{k=1}^{\nu_m}\; \lambda_k\, \operatorname{id}_{V(\chi_k )}
\; \in \,\bigoplus_{k=1}^{\nu_m}\; \operatorname{End}(V (\chi_k ))
\;(\subset\, \operatorname{End}(V_m)).
$$
In terms of the identification (3.1), we can write the Hermitian metric
$(\;,\;)_m$ on $\frak h_m$ in (2.1)  in the form
$$
(\lambda, \mu )_m:= \Sigma_{k=1}^{\nu_m} \, n_k \,\lambda_k \bar{\mu}_k/m^{n+2},
$$
where $\lambda = (\lambda_1, \lambda_2,\dots, \lambda_{\nu_m})$ and 
$\mu = (\mu_1, \mu_2,\dots, \mu_{\nu_m})$ are in $\Bbb C^{\nu_m}$.
By the identification (3.1), corresponding to the decomposition
$\frak h_m = \frak t \oplus \frak t^{\perp}$, we have the 
orthogonal direct sum
$$
\Lambda_m \; =\; \Lambda (\frak t ) \oplus \Lambda (\frak t^{\perp}),
$$
where $\Lambda  (\frak t )$ and $\Lambda  (\frak t^{\perp})$ are the subspace of 
$\Lambda_m$ associated to $\frak t$ and $\frak t^{\perp}$, respectively.
Take a Hermitian metric $\rho_k$ on $V (\chi_k )$, and 
for the metric  
$$
\rho \; :=\; \bigoplus_{k=1}^{\nu_m}\; \rho_k
$$ 
on $V_m$,  we see that
$V(\chi_k ) \perp V(\chi_{k'} ) $ whenever $k \neq k'$. 
By choosing an orthonormal basis $\{\,s_{k,i}\,;\, i=1,2,\dots, n_k\,\}$ for the Hermitian vector space $(V(\chi_k ), \rho_k )$, we now set
$$
j (k, i) \; :=\; i\,+ \sum_{l =1}^{k-1} \; n^{}_{l},
\;\; i=1,2,\dots, n_k;\; k=1,2,\dots, \nu_m,
\leqno{(3.2)}
$$
where the right-hand side denotes $i$ in the special case $k =1$.
By writing $s_{k,i}$ as $s_{j (k,i)}$, we have an orthonormal basis 
$$
\mathcal{S} := \{s_1, s_2, \dots , s_{N_m}\}
$$ 
for  $(V_m , \rho )$.  By this basis, the vector space $V_m$ and
the algebraic group $\operatorname{SL}(V_m)$ are identified with 
$\Bbb C^{N_m} = \{ \,(z_1,\dots,z_{N_m} )\,\}$ and 
$\operatorname{SL}(N_m,\Bbb C )$, respectively.
In terms of $\mathcal{S}$, the Kodaira embedding $\Phi_m$ is given
by 
$$
\Phi_m (x) := (s_1(x):\dots:s_{N_m}(x)), \qquad x\in M.
$$
Consider the 
associated {\it Chow norm} $W^*_m \owns \xi \mapsto \|\xi\|^{}_{\operatorname{CH}(\rho )}\in \Bbb R_{\geq 0}$ 
as in  Zhang \cite{Zh} 
(see also \cite{M0}).
Then by the closedness of $G_m \cdot \hat{M}_m$ in $W_m^*$ (cf. (2.2) and Definition 2.3), 
the Chow norm on the orbit 
$G_m \cdot \hat{M}_m$ takes its minimum at $g_m\cdot \hat{M}_m$ for some  $g_m\in G_m$.  
Note that,  by complexifying
$$
K_m\, :=\,\prod_{k=1}^{\nu_m}\,\operatorname{SU}(V(\chi_k )),
$$
we obtain the reductive algebraic group $S_m$.
For each $\kappa\in K_m$
and each diagonal matrix $\Delta$ in $\frak{sl}(N_m, \Bbb C )$,
we put
$$
e(\kappa, \Delta ) := \exp\{\operatorname{Ad}
(\kappa ) \Delta \}.
$$ 
Then $g_m$ is written as $\kappa_1 \cdot  e(\kappa_0,  {D} )
$ for some $\kappa_0$, $\kappa_1 \in K_m$
and a 
diagonal matrix ${D} = ({d}_j )_{1\leq j\leq N_m}$ 
 in $\frak{sl}(N_m, \Bbb C )$ with the $j$-th diagonal element ${d}_j$. 
Put ${g}'_m := e(\kappa_0, {D})$.
In view of
$\| {g}'_m\cdot \hat{M}_m \|^{}_{\operatorname{CH}(\rho )} = 
\| g_m \cdot \hat{M}_m \|^{}_{\operatorname{CH}(\rho )}$ (cf. \cite{Sa}, Prop.~4.1), 
we obtain
$$
\|  {g}'_m\cdot \hat{M}_m \|^{}_{\operatorname{CH}(\rho )}
\leq \| e(\kappa_0, t (X^{}_{\lambda}+A))\cdot 
{g}'_m\cdot \hat{M}_m \|^{}_{\operatorname{CH}(\rho )}, \;\;\, t \in \Bbb C,
\leqno{(3.3)}
$$
for all $\lambda = (\lambda_1, \lambda_2, \dots, \lambda_{\nu_m}) \in \Lambda (\frak t^{\perp})$ and all diagonal matrices $A = (a_j)_{1\leq j \leq N_m}$ in $\frak s_m$,
where each $a_j$ denotes the $j$-th diagonal element of $A$.
We now write $a_{j(k,i)}$ as $a_{k,i}$ for simplicity. 
Put
$$
{s}'_j:= \; \kappa_0^{-1}\cdot s_j,\quad
b_{k,i}:= \; \lambda_k + a_{k,i}, \quad
c_{k,i} := \; \exp {d}_{j (k,i)}.
$$
Then we shall now identify $V_m$ with 
$\Bbb C^{N_m} = \{\,(z'_1, \dots, z'_{N_m})\,\}$ by the orthonormal basis
${\mathcal{S}}':=\{{s}'_1, {s}'_2, \dots, {s}'_{N_m}\}$ for $V_m$.
In view (3.2), we rewrite $s'_j$, $z'_j$ as 
$s'_{k,i}$, $z'_{k,i}$, respectively by
$$
{s}'_{k,i}:= {s}'_{j (k,i)},
\quad
z'_{k,i}:= z'_{j (k,i)},
$$
where $k = 1,2,\dots,\nu_m$ and $i = 1,2,\dots,n_k$.
By writing $b_{k,i}$, $c_{k,i}$ also as $b_{j (k,i)}$, $c_{j (k,i)}$, 
respectively,
we consider the diagonal matrices $B$ and $C$ of order $N_m$
with the $j$-th diagonal elements $b_j$ and $c_j$, respectively.
Note that the right-hand side of (3.3) is 
$$
\| (\exp tB) \cdot C \cdot
\kappa_0^{-1}\cdot \hat{M}_m \|^{}_{\operatorname{CH}(\rho )},
$$
and its derivative at $t =0$ vanishes 
by virtue of the inequality (3.3). Hence, by setting $\Theta := (\sqrt{-1}/2\pi )\partial\bar{\partial}\log 
(\Sigma_{k=1}^{\nu_m} \Sigma_{i=1}^{n_k}|c_{k,i}z'_{k,i}|^2)$, 
we obtain the equality (see for instance  (4.4) in \cite{M0})
$$
\int_M 
\frac{\Sigma_{k=1}^{\nu_m} \Sigma_{i=1}^{n_k}b_{k,i} |c_{k,i}{s}'_{k,i}|^2}
{\Sigma_{k=1}^{\nu_m} \Sigma_{i=1}^{n_k} |c_{k,i}{s}'_{k,i}|^2}
\,{\Phi}_m'^*(\Theta^n ) \; =\; 0
\leqno{(3.4)}
$$
for all $\lambda \in \Lambda (\frak t^{\perp} )$ 
and all diagonal matrices $A$ in the Lie algebra $\frak s_m$, 
where ${\Phi}'_m : M \hookrightarrow \Bbb P^*(V_m)$ is the Kodaira embedding 
of $M$ by ${\mathcal{S}}'$ which
sends each $x\in M$ to 
$({s}'_1(x):{s}'_2(x): \dots :{s}'_{N_m}(x))\in \Bbb P^*(V_m)$. Here we regard
$$
s'_{k,i}\; =\; {\Phi'}^*z'_{k,i}.
$$
Let ${k}_0\in \{1,2,\dots,\nu_m\}$ and let $i_1,i_2 \in \{1,2,\dots,n^{}_{{k}_0}\}$ with 
$i_1 \neq i_2$. Using Kronecker's delta, we first specify the real diagonal matrix $B$ 
by 
$$
\lambda_k = 0
\quad\text{ and }\quad a_{k,i} = \delta_{k {k}_0} (\delta_{i i_1} - \delta_{i i_2} ),
$$
where $k=1,2,\dots,\nu_m; \,i=1,2,\dots,n_k$. By (3.4) applied to this $B$, and let $(i_1,i_2)$ run through the set of all pairs of two distinct integers in $\{1,2,\dots,n_{{k}_0}\}$, where positive integer ${k}_0$ varies from $1$ to $\nu_m$.
Then there exists a positive constant $\beta_k > 0$ independent of the choice of 
$i$ in $\{1,2,\dots,n_k\}$ such that, for all $i$,
$$
\int_M \frac{|c_{k,i}{s}'_{k,i}|^2}
{\Sigma_{k=1}^{\nu_m} \Sigma_{i=1}^{n_k}|c_{k,i}{s}'_{k,i}|^2}
\,{\Phi}_m'^*(\Theta^n ) \; =\; \beta_k,
\;\; k=1,2,\dots,\nu_m.
\leqno{(3.5)}
$$
Let ${k}_0, i_1, i_2$ be as above, and let $\kappa_2$ be the element in $K_m$ such that 
$\kappa_2 \,z'_{{k}_0,i_1} = (1/\sqrt{2})(z'_{{k}_0,i_1} - z'_{{k}_0,i_2})$, 
$\kappa_2 \,z'_{{k}_0,i_2} = (1/\sqrt{2})(z'_{{k}_0,i_1} + z'_{{k}_0,i_2})$
and that $\kappa_2$ fixes all other $z_{k,i}$'s. 
Let $\kappa_3$ be the element in $K_m$ such that 
$\kappa_3 \,z'_{{k}_0,i_1} =(1/\sqrt{2})(z'_{{k}_0,i_1} + \sqrt{-1}\, z'_{{k}_0,i_2})$, 
$\kappa_3 \,z'_{{k}_0,i_2} = (1/\sqrt{2})(\sqrt{-1}\,z'_{{k}_0,i_1}
  + z'_{{k}_0,i_2})$
and that $\kappa_3$ fixes all other $z'_{k,i}$'s.  Now 
$$
\| \kappa_2 {g}'_m \cdot \hat{M}_m \|_{\operatorname{CH}(\rho )} 
= \| \kappa_3 {g}'_m \cdot \hat{M}_m \|_{\operatorname{CH}(\rho )} 
=\|  {g}'_m \cdot \hat{M}_m \|_{\operatorname{CH}(\rho )} ,
$$
and note that 
$$
2z'_{{k}_0,i_1}\bar{z}'_{{k}_0,i_2} = (|\kappa_2 \,z'_{{k}_0,i_2}|^2 - |\kappa_2 \,z'_{{k}_0,i_1}|^2)- \sqrt{-1}\,(|\kappa_3 \,z'_{{k}_0,i_2}|^2 - |\kappa_3 \,z'_{{k}_0,i_1}|^2).
$$
Hence replacing ${g}'_m $ by 
$\kappa_{\alpha} {g}'_m $, $\alpha =2,3$, 
in (3.3),
we obtain the case $k'=k''$ of the following by an argument as in deriving (3.5) from (3.3):
$$
\int_M \frac{{s}'_{k',i'}\bar{s}'_{k'',i''}}
{\Sigma_{k=1}^{\nu_m} \Sigma_{i=1}^{n_k}|c_{k,i}{s}'_{k,i}|^2}
\,{\Phi}_m'^*(\Theta^n ) \; =\; 0,
\;\;\text{ if $(k',i') \neq (k'',i'')$}.
\leqno{(3.6)}
$$
Here (3.6) holds easily for $k'\neq k''$, since for every element $g$ 
of the maximal compact subgroup of $T$, we have:
\begin{align*}
\text{L.H.S. of (3.6)}\; &=\;\int_M \,
g^*\left \{ \frac{{s}'_{k',i'}\bar{s}'_{k'',i''}}
{\Sigma_{k=1}^{\nu_m} \Sigma_{i=1}^{n_k}|c_{k,i}{s}'_{k,i}|^2}
\,{\Phi}_m'^*(\Theta^n )\right \}\\
&= \;\frac{\chi_{k''}(g)}{\chi_{k'}(g)}\,
\int_M 
\frac{{s}'_{k',i'}\bar{s}'_{k'',i''}}
{\Sigma_{k=1}^{\nu_m} \Sigma_{i=1}^{n_k}|c_{k,i}{s}'_{k,i}|^2}
\,{\Phi}_m'^*(\Theta^n ).
\end{align*}
Put $\beta := (\beta_1, \beta_2,\dots,\beta_{\nu_m})\in \Bbb R^{\nu_m}$ 
and $\beta_0:= (\Sigma_{k=1}^{\nu_m} n_k \beta_k )/N_m$, where $\beta_k$ is given in (3.5).   In view of 
$N_m = \Sigma_{k=1}^{\nu_m} n_k$, by setting $\underline{\beta}_k:= \beta_k -\beta_0$,
we have
$$
\underline{\beta} := (\underline{\beta}_1, \underline{\beta}_2, \dots, \underline{\beta}_{\nu_m})
\in \Lambda_m.
$$
Next for each $\lambda\in \Lambda (\frak t^{\perp} )$, by setting $a_{k,i} = 0$ for all $(k,i)$, 
the equality (3.4) above implies
$\,0\, =\, \Sigma_{k=1}^{\nu_m} \, (n_k \lambda_k)\,\beta_k $.
From this  together with the equality $\Sigma_{k=1}^{\nu_m} \,n_k \lambda_k =0$, we obtain 
$(\lambda, \underline{\beta} )_m =0$, i.e.,
$$
\underline{\beta} \in \Lambda (\frak t ).
\leqno{(3.7)}
$$
We now define a Hermitian metric $h_{\operatorname{FS}}$ (cf. \cite{Zh})
for $L^{m}$  as follows. Let $u$ be a local section for $L^m$. Then${}^{\dagger}$
$$
|u|^2_{h_{\operatorname{FS}}}\, := \;
\frac{|u|^2}
{\Sigma_{k=1}^{\nu_m} \Sigma_{i=1}^{n_k}\; |c_{k,i}{s}'_{k,i}|^2}.
\leqno{(3.8)}
$$
\footnotetext{ ${}^{\dagger}$In view of (3.8), there is some error in \cite{M0}. 
Actually, for the numerator of (5.7) in the paper \cite{M0}, please read $(N_m + 1) |s|^2$.}
For the Hermitian metric $h_m := (h_{\operatorname{FS}})^{1/m}$ for $L$, 
we consider the associated K\"ahler metric $\omega_m := c_1(L; h_m)_{\Bbb R}$
 on $M$.
In view of (3.5),
$$
\begin{cases}
\;\;\,\beta_0 & =\, (\Sigma_{k=1}^{\nu_m}n_k\beta_k)/N_m\,
\,=\, N_m^{-1}\int_M\Phi_m^*(\Theta^n )\\
&=\,  N_m^{-1}m^n c_1(L)^n[M]\, = \, n!\, \{1+ O(1/m)\}.
\end{cases}
\leqno{(3.9)}
$$
Then for 
$\gamma_{m,k}:= \beta_k/\beta_0$ and 
$\sigma_{k,i}:= c_{k,i}\,{s}'_{k,i}(m^n\beta^{-1}_k)^{1/2}$, we have
$$
\begin{cases}
&\Sigma_{k=1}^{\nu_m}\Sigma_{i=1}^{n_k}\, \gamma_{m,k} |\sigma_{k,i}|_{h_m}^2\,
=\,\Sigma_{k=1}^{\nu_m}\Sigma_{i=1}^{n_k}\, \gamma_{m,k}|\sigma_{k,i}|_{h_{\operatorname{FS}}}^2\\
&=\, (m^n/\beta_0 )\,\Sigma_{k=1}^{\nu_m}\Sigma_{i=1}^{n_k}\, |c_{k,i} {s}'_{k,i}|_{h_{\operatorname{FS}}}^2\, =\, m^n/\beta_0.
\end{cases}
\leqno{(3.10)}
$$
By operating $(\sqrt{-1}/2\pi ) \bar{\partial}{\partial}\log$ on (3.8), we obtain 
$$
{\Phi}_m'^*(\Theta ) \;=\; c_1 (L^{m}; h_{\operatorname{FS}})_{\Bbb R} \; =\; m\, c_1(L;h_m)_{\Bbb R}\; = \;  m\, \omega_m.
\leqno{(3.11)}
$$
Then in terms of the Hermitian $L^2$ inner product (1.1), we see from (3.5), (3.6), (3.8) and (3.11) that $\{\sigma_{k,i}\,;\, k=1,\dots,\nu_m,\, i=1,\dots, n_k\}$ is an orthonormal basis for $V_m$. 
Moreover, (3.10) is rewritten as 
$$
B^{\circ}_m(\omega_m) \; =\; m^n/\beta_0, 
\leqno{(3.12)}
$$ 
where $B^{\circ}_m(\omega_m)$ is as in (1.2). Hence $\omega_m$ is a polybalanced metric,
and the proof of Theorem A is reduced to showing (1.3) and (1.4). By summing up, we obtain

\medskip\noindent
{\bf Theorem C:} 
{\em 
If $(M, L^m)$ is Chow-stable relative to 
$T$ for a positive integer $m$, then the K\"ahler class $c_1(L)_{\Bbb R}$ admits a polybalanced metric $\omega_m$ 
with the weights $\gamma_{m,k}$ as above.
}

\section{The asymptotic behavior of the weights $\gamma_{m,k}$}

The purpose of this section is to prove (1.3). If $\underline{\beta} = 0$, then we are done.
Hence, we may assume that $\underline{\beta} \neq  0$.
Consider the sphere 
$$
\Sigma\; :=\;\{\,X\in \frak t_{\Bbb R}\,; \,\langle X,X\rangle_0 =1\,\}
$$ 
in $\frak t_{\Bbb R}:= \frak t \cap (\frak h_m)_{\Bbb R}$, 
where $\langle \;,\,\rangle_0$ denotes the positive definite symmetric bilinear form 
on $\frak g$ as in \cite{FM}.
 Since all components $\underline{\beta}_k$ of 
$\underline{\beta}$ are real, 
we see from (3.7) that, in view of (3.1), 
$\lambda:= r_m\,\underline{\beta}$ satisfies
$$
X_{\lambda} \;\in\; \Sigma 
\leqno{(4.1)}
$$
for some positive real number $r_m$. 
Hence by writing $\lambda = (\lambda_1,\lambda_2,\dots,\lambda_{\nu_m})$, 
we obtain positive constants $C_2$, $C_3$ independent of $k$ and $m$ such that
(see for instance \cite{Ma}, Lemma 2.6)
$$
-\,C_2 \,m \; \leq \;  \lambda_k\; \leq \; C_3\, m.
\leqno{(4.2)}
$$
Put $g(t):= \exp (tX_{\lambda})$ and $\gamma (t) := \log \| g(t)\cdot C\cdot \kappa_0^{-1}\cdot \hat{M}_m
\|_{\operatorname{CH}(\rho )}$, $t \in \Bbb R$, 
by using the notation in Sect.~3. Since $g(t)$ commutes with $C\cdot \kappa_0^{-1}$,
we see that $g(t)$ defines a holomorphic automorphism of $\Phi_m'(M)$.
In view of Theorem 4.5 in \cite{M0}, it follows from Remark 4.6 in  \cite{M0} that
(cf. \cite{Zh})
$$
\dot{\gamma} (t)\; =\; 
\dot{\gamma} (0) \; =\; \Sigma_{k=1}^{\nu_m}\, n_k \lambda_k \beta_k ,
\qquad -\infty < t< +\infty.
\leqno{(4.3)}
$$
Consider the classical Futaki invariant $\mathcal{F}_1 (X_{\lambda})$ associated to the holomorphic vector 
field $X_{\lambda}$ on $(M,L)$.
Since $M$ is smooth, this coincides with the 
corresponding Donaldson-Futaki's invariant for test configurations. 
Then by applying
Lemma 4.8 in \cite{M} to the product configuration of $(M,L)$ associated to the 
one-parameter group generated by $X_{\lambda}$ on the central fiber, 
we obtain (see also \cite{D2}, \cite{Sz})
$$
\lim_{t\to -\infty} \dot{\gamma} (t)\; =\;
\,C_4\,\{ \mathcal{F}_1(X_{\lambda}) + O(1/m)\}\, m^n,
\leqno{(4.4)}
$$
where $C_4 :=\, (n+1)!\,c_1(L)^n[M] > 0$.
Hence by %$\underline{\beta}_k = {\beta}_k - \beta_0$ and 
$\Sigma_{k=1}^{\nu_m} \,n_k \lambda_k  
%=r_m\,\Sigma_{k=1}^{\nu_m}\, n_k \underline{\beta}_k 
= 0$ and $\lambda =\, r_m\, \underline{\beta}$, it now follows from (4.3) and (4.4) that
$$
\begin{cases}
\quad C_4\,\{ \mathcal{F}_1(X_{\lambda}) + O(1/m)\}\, m^n 
&= \; r_m^{-1}\,\Sigma_{k=1}^{\nu_m}\, n_k\lambda_k^{2}\\
&=\; r_m^{-1}\,m^{n+2}\,\langle X_{\lambda}, X_{\lambda} \rangle_m,
\end{cases}
\leqno{(4.5)}
$$
where by (4.1) above, (7) in \cite{D2} (see also \cite{Sz}) implies 
$\langle X_{\lambda}, X_{\lambda} \rangle_m \geq C_5$ 
for some positive real constant $C_5$ independent of $m$.
Furthermore $\mathcal{F}_1(X_{\lambda})  = O(1)$ 
again by (4.1).
Hence from (4.5), we obtain
$$
r_m^{-1} \; = \; O(1/m^2).
\leqno{(4.6)}
$$
In view of (3.9), since $\underline{\beta}_k = r_m^{-1}\lambda_k$, 
(4.2) and (4.6) imply the required estimate (1.3) as follows:
$$
\gamma_{m,k} -1 \; =\; (\beta_k/\beta_0) \,-\,1 \;=\; \underline{\beta}_k/\beta_0
\; = \; O(1/m).
$$

\section{Proof of (1.4) and Corollary B}

In this section, keeping the same notation as in the preceding sections,
we shall prove (1.4) and Corollary B. 

\medskip\noindent
{\em Proof of $(1.4)$}:  
For $X_{\lambda}$ in (4.1), the associated Hamiltonian function $f_{\lambda} 
\in C^{\infty}(M)_{\Bbb R}$ on the K\"ahler manifold $(M, \omega_m )$ is
$$
f_{\lambda} = \frac{\Sigma_{k=1}^{\nu_m} \Sigma_{i=1}^{n_k}\lambda_k
\gamma_{m,k} |\sigma_{k,i}|_{h_m}^2}
{m\Sigma_{k=1}^{\nu_m} \Sigma_{i=1}^{n_k} \gamma_{m,k} |\sigma_{k,i}|_{h_m}^2}
= \frac{\Sigma_{k=1}^{\nu_m} \Sigma_{i=1}^{n_k}\lambda_k
|c_{k,i}s'_{k,i}|^2}
{m\Sigma_{k=1}^{\nu_m} \Sigma_{i=1}^{n_k} 
 |c_{k,i}s'_{k,i}|^2},
$$
where by (4.2), when $m$ runs through the set of all sufficiently large integers, the function $f_{\lambda}$
is uniformly $C^0$-bounded. Now by (3.10),
$$
\Sigma_{k=1}^{\nu_m} \Sigma_{i=1}^{n_k}\,\lambda_k
\gamma_{m,k} |\sigma_{k,i}|_{h_m}^2 \; =\; m^{n+1}f_{\lambda}/\beta_0.
\leqno{(5.1)}
$$
We now define a function $I_m$ on $M$ by
$$
I_m\,  :=\;(r_m^{-1}/\beta_0)\,\Sigma_{k=1}^{\nu_m} \Sigma_{i=1}^{n_k}\, 
\{1 - (\gamma_{m,k})^{-1}\}\, \lambda_k
 \gamma_{m,k}\, |\sigma_{k,i}|_{h_m}^2.
 \leqno{(5.2)}
$$
Then by (1.3), (3.9), (3.10), (4.2) and (4.6), 
we easily see that 
$$
I_m\, =\;O(m^{-2}\Sigma_{k=1}^{\nu_m} \Sigma_{i=1}^{n_k} \gamma_{m,k}\, |\sigma_{k,i}|_{h_m}^2)\; =\; O(m^{n-2}).
\leqno{(5.3)}
$$
By (3.10) together with (5.1) and (5.2), it now follows that
\begin{align*}
\Sigma_{k=1}^{\nu_m} \Sigma_{i=1}^{n_k}
 |\sigma_{k,i}|_{h_m}^2  &=\Sigma_{k=1}^{\nu_m} \Sigma_{i=1}^{n_k}\gamma_{m,k}
 |\sigma_{k,i}|_{h_m}^2 
  -\Sigma_{k=1}^{\nu_m} \Sigma_{i=1}^{n_k}(\gamma_{m,k}-1)
 |\sigma_{k,i}|_{h_m}^2\\
 &= (m^n/\beta_0) \, -\, \Sigma_{k=1}^{\nu_m} \Sigma_{i=1}^{n_k}\, (\underline{\beta}_k/\beta_0)\,
 |\sigma_{k,i}|_{h_m}^2\\
 &= (m^n/\beta_0) \,+ \,  I_m \,
  -\, (r_m^{-1}/\beta_0)\,\Sigma_{k=1}^{\nu_m} \Sigma_{i=1}^{n_k}\,  \lambda_k
 \gamma_{m,k}\, |\sigma_{k,i}|_{h_m}^2\\
 &= (m^n/\beta_0) \,+ \,  I_m \,
  -\, (r_m^{-1}m^2/\beta^{\,2}_0)\,f_{\lambda}\,m^{n-1}.
\end{align*}
By (3.9), $m^n/\beta_0 = N'_m$. 
Moreover, by (4.6), $r_m^{-1}m^2/\beta^{\,2}_0 = O(1)$. Since
$$
f_m \; :=\; - \,  (r_m^{-1}m^2/\beta^{\,2}_0)\,f_{\lambda}, \quad m\gg 1,
$$
are uniformly $C^0$-bounded Hamiltonian functions on 
$(M, \omega_m )$ associated to holomorphic 
vector fields in $\frak t$, in view of (5.3),  we obtain
$$
B_m^{\bullet}(\omega_m ) \; =\; \Sigma_{k=1}^{\nu_m} \Sigma_{i=1}^{n_k}
 |\sigma_{k,i}|_{h_m}^2 \; =\; N'_m \, +\,f_m \, m^{n-1}\, + \, O(m^{n-2}),
$$
as required.

\medskip\noindent
{\em Proof of Corollary B}:  Since the classical Futaki character $\mathcal{F}_1$ vanishes 
on $\frak t$,
we have $\mathcal{F}_1(X_{\lambda}) = 0$ in (4.5), so that
$$
r_m^{-1} \; =\; O(1/m^3).
$$
Then from (3.9) and  $\underline{\beta}_k =\, r_m^{-1}\,\lambda_k$,
by looking at (4.2), we obtain the following required estimate:
$$
\gamma_{m,k} -1 \; =\; \underline{\beta}_k/\beta_0
\; = \; O(1/m^2).
$$
Hence $B^{\bullet}_m (\omega_m )\, =\, \{1+ O(1/m^2)\}\,B^{\circ}_m(\omega_m )$. Integrating this over $M$ by the volume form $\omega_m^n$, in view of (3.12), we see that
$$
N'_m \; =\; \{1+ O(1/m^2)\} (m^n/\beta_0 ).
$$
Therefore, from (3.12) and $B^{\bullet}_m (\omega_m )\, =\, \{1+ O(1/m^2)\}\,B^{\circ}_m(\omega_m )$, we now conclude that $B^{\bullet}_m (\omega_m )\, =\, N'_m + O (m^{n-2})$, as required.

%%%%%%%%%%%%%%%%%%%%%%%%%%%%%
\bigskip\noindent
{\footnotesize
{\sc Department of Mathematics}\newline
{\sc Osaka University} \newline
{\sc Toyonaka, Osaka, 560-0043}\newline
{\sc Japan}}
%%%%%%%%%%%%%%%%%%%%%%%%%%%%%
\end{document}